\documentclass{article}[12pt]

\def\proof{{\sc Proof. }}
\usepackage{amsfonts}
\usepackage{amssymb}

\newtheorem{thm}{Theorem}
\newtheorem{prop}[thm]{Proposition}
\newtheorem{lemma}[thm]{Lemma}

\def\R{\mathbb{R}}
\def\H{\mathbb{H}}
\def\Re{{\rm Re}}
\def\Im{{\rm Im}}
\def\HP{{\H}{\rm P}}
\def\RP{\R{\rm P}}
\def\d{\partial}
\def\<{\langle}
\def\>{\rangle}
\def\Box{\square}

\title{Maps That Take Lines to Circles, in Dimension 4}
\author{V. Timorin}
\date{}

\begin{document}
\maketitle

{\small {\bf Abstract.} 
We list all analytic diffeomorphisms between an open subset of the 4-di\-men\-sional 
projective space and an open subset of the 4-dimensional sphere that take all 
line segments to arcs of round circles. 
These are the following: restrictions of the quaternionic Hopf fibrations and 
projections from a hyperplane to a sphere from some point. 
We prove this by finding the exact solutions of the corresponding system of 
partial differential equations.}

\section{Introduction}

Let $U$ be an open subset of the 4-dimensional real projective space $\RP^4$
and $V$ an open subset of the 4-dimensional sphere $S^4$. We study diffeomorphisms
$f:U\to V$ that take all line segments lying in $U$ to arcs of round circles lying
in $V$. For the sake of brevity we will always say in the sequel that $f$ {\em takes
all lines to circles}. The purpose of this article is to give the complete list of 
such analytic diffeomorphisms. 

\paragraph{Remark.}
Given a diffeomorphism $f:U\to V$ that takes lines to circles, we can compose 
it with a projective transformation in the preimage (which takes lines to lines)
and a conformal transformation in the image (which takes circles to circles).
The result will be another diffeomorphism taking lines to circles. 

\paragraph{Example 1.}
For example, suppose that $S^4$ is embedded in $\R^5$ as a Euclidean sphere and
take an arbitrary hyperplane and an arbitrary point in $\R^5$. Obviously, the projection
of the hyperplane to $S^4$ form the given point takes all lines to circles. 
The restriction of this projection to some open subset of the hyperplane is a
diffeomorphism between this subset and its image. Diffeomorphisms obtained in this way 
will be called {\em classical projections}. Of course, classical projections 
live in all dimensions, not only 4.

\paragraph{Example 2.}
Another example comes from the {\em quaternionic Hopf fibration}. Let us recall the 
definition. Consider the standard projection from the left (resp., right) quaternionic
2-dimensional vector space $\H^2$ to the left (resp., right) quaternionic projective
line $\HP^1$. The latter is identified with $S^4$, a conformal quaternionic coordinate being the 
ratio of the homogeneous coordinates (in a specified order). Clearly, this projection 
descends to $\RP^7$ --- the real projectivization of $\H^2$. Thus we obtain a map 
from $\RP^7$ to $S^4$. This is the quaternionic Hopf fibration. It is known 
(see e.g. \cite{L}) that 
the quaternionic Hopf fibration takes all lines to circles. Therefore, the same is 
true for the restriction of it to any 4-dimensional projective subspace of $\RP^7$.

Our main result is the following

\begin{thm}
\label{main}
Suppose that an analytic diffeomorphism $f$ between an open set in $\RP^4$ and an open set in
$S^4$ takes all lines to circles. Then it 
is either a restriction of a classical projection or a restriction of a (left or right)
quaternionic Hopf fibration.
\end{thm}

{\sc Remark.} When we say that a map $f:U\subseteq\RP^4\to S^4$ is a restriction of 
a classical projection, we mean the following. There is a projective map $i$ from 
$\RP^4$ to $\R^5$ defined everywhere on $U$, a conformal identification $j$ of $S^4$ 
with a Euclidean sphere in $\R^5$ and a central projection $\pi$ to this Euclidean sphere 
from some point such that $f=j^{-1}\circ\pi\circ i$. 
Analogously, the statement ``$f$ is a restriction of a quaternionic Hopf fibration''
has the following meaning. There is a projective map $i$ form $\RP^4$ to $\RP^7$,
a (right or left) quaternionic Hopf fibration $\pi:\RP^7\to S^4$ and a conformal 
identification $j$ of the sphere in the image of $f$ with the sphere in the image 
of $\pi$ such that $f=j^{-1}\circ\pi\circ i$.
Note in particular, that Theorem \ref{main} is insensitive to a projective 
transformation in the preimage and a conformal transformation in the image.  

Actually, it is not hard to see that we can always assume $j$ to be fixed from 
the very beginning.

\paragraph{History of the problem.}
The problem of finding maps that take lines to circles came from nomography 
(for an introduction to nomography see \cite{GKh}). G.S. Khovanskii in 1970-s 
posed the following question: find all diffeomorphisms between open subsets of 
$\R^2$ that take lines to circles, or, in the language of nomography, that 
transform nomograms with aligned points to circular nomograms. Circular nomograms 
are more convenient in practice whereas nomograms with aligned points are 
theoretically easier to deal with. This problem was solved by A.G. Khovanskii 
\cite{Kh} who proved that all such diffeomorphisms come from projections 
of planes in $\R^3$ to Euclidean spheres. Actually, his result was stated in 
a different form: up to a projective transformation in the preimage and a 
M\"obius transformation in the image there are only 3 such diffeomorphisms, 
and they correspond to classical geometries. Izadi in \cite{Iz} extended the 
results of Khovanskii to the 3-dimensional case. 

It turned out (see \cite{Tim1}) that Khovanskii's theorem does not extend to 
dimension 4. The simplest counterexample is a complex projective transformation
which takes lines to circles but does not come from a projection of a hyperplane
to a sphere. This is a particular case of example 2 above. In \cite{Tim1} 
it was proved that all ample enough rectifiable bundles of circles passing 
through a point in $\R^4$ are obtained by means of example 2 (those  obtained
from example 1 can be obtained from example 2 as well). This 
may be regarded as a first step in proving Theorem \ref{main}. This article 
completes the proof. 

\paragraph{Outline of the proof of Theorem \ref{main}.} 
Suppose that $f:U\to V$ is a diffeomorphism
taking all lines to circles. First choose a projective identification of $U$ with 
a region in $\R^4$ and a conformal identification of $V$ with a region in $\H$,
the skew-field of quaternions. Then $f$ can be regarded as a quaternion-valued function
on $U$.

For any constant vector field $\alpha$ denote by $\d_\alpha$ the Lie derivative 
along $\alpha$. Put $A_\alpha=\d_\alpha f$. The quaternions $A_\alpha$ are 
best to be understood as components of a quaternion-valued differential 1-form $A$
on $U$, namely, the differential of $f$. This form is closed: $dA=0$ or, in 
components, $\d_\alpha A_\beta=\d_\beta A_\alpha$ for any pair of constant 
vector fields $\alpha$ and $\beta$ on $U$. 

By Lemma 5.2 from \cite{Tim1}, at every point of $U$ there exists a quaternion 
$B_\alpha$ that depends linearly on $\alpha$ and satisfies one of the equations
$$\d_\alpha A_\alpha=B_\alpha A_\alpha\quad{\rm or}\quad \d_\alpha A_\alpha=A_\alpha B_\alpha.$$
The products in the right-hand sides are in the sense of quaternions.
Assume that the first equation holds everywhere on $U$ (this is not very
restrictive: see the end of Section 4).
We show in Section 2 that the integrability condition for this equation is
$$\d_\alpha B_\beta=\frac 12 B_\alpha B_\beta+\frac 13 C(\overline{A_\alpha}A_\beta+
\overline{A_\beta}A_\alpha+A_\alpha\overline{A_\beta}).$$
where $C$ is a quaternion-valued function. 
We call it the {\em first integrability condition}. Then we use once again that $f$
takes lines to circles, now to deduce that $C$ is real-valued. If $C=0$ identically,
then the first integrability condition is integrable, and it gives us the maps from 
example 2. This will be proved in Section 3.

Suppose now that $C\ne 0$. Then we need to write down an integrability condition for 
the first integrability condition --- the {\em second integrability condition}.
It will be obtained in Section 4.
The second integrability condition turns out to be a differential relation between
$A$, $B$ and $C$ only, with no additional parameters. It expresses all first logarithmic 
derivatives of $C$ in 
terms of the values of $A$ and $B$. The fact that $C$ is real imposes a restriction on 
possible values of $A$ and $B$. Namely, $B_\alpha=p(A_\alpha)+A_\alpha q$
where $p$ is a real-valued 1-form and $q$ is a quaternion. Now fix 
some admissible values of $A$, $B$ and $C$ at some point of $U$. We will show
in Section 5
that there exists a map from example 1 with the same values of $A$, $B$
and $C$ at the given point. On the other hand, from the second integrability 
condition it follows that such analytic map is unique. Thus our map $f$ is from 
example 1.

\section{First integrability condition}

In this section, we obtain an integrability condition for the equation
$$\d_\alpha A_\alpha=B_\alpha A_\alpha.\eqno{(1)}$$

Since $A_\alpha$ are derivatives, we have $\d_\alpha A_\beta=\d_\beta A_\alpha$. 
From this and from equation (1) it follows immediately that 
$$\d_\alpha A_\beta=\frac 12(B_\alpha A_\beta+B_\beta A_\alpha).\eqno{(2)}$$
 
An integrability condition for equation (1) is obtained from the 
equality $\d_\alpha\d_\beta A_\gamma=\d_\beta\d_\alpha A_\gamma$ by expressing derivatives 
of $A$ with the help of (2). If we introduce a tensor $C$ given in components by
$$C_{\alpha\beta}=\d_\alpha B_\beta-\frac 12 B_\alpha B_\beta,$$
then the integrability condition reads as follows:
$$C_{\alpha\beta}A_\gamma+C_{\alpha\gamma}A_\beta=C_{\beta\alpha}A_\gamma+
C_{\beta\gamma}A_\alpha.\eqno{(3)}$$
If we apply to equation (3) the transposition of indices $\beta\leftrightarrow\gamma$,
then we have: 
$$C_{\alpha\beta}A_\gamma+C_{\alpha\gamma}A_\beta=C_{\gamma\alpha}A_\beta+C_{\gamma\beta}A_\alpha.
\eqno{(3')}$$
Take the sum of equations (3) and $(3')$. It can be written as follows: 
$$(C_{\beta\gamma}+C_{\gamma\beta})A_\alpha=(2C_{\alpha\gamma}-C_{\gamma\alpha})A_\beta+
(2C_{\alpha\beta}-C_{\beta\alpha})A_\gamma.\eqno{(3'')}$$

Put $\beta=\gamma$ in equation (3). We obtain that
$$(2C_{\alpha\beta}-C_{\beta\alpha})A_\beta=C_{\beta\beta}A_\alpha.\eqno{(4)}$$
or, replacing $\beta$ with $\gamma$, 
$$(2C_{\alpha\gamma}-C_{\gamma\alpha})A_\gamma=C_{\gamma\gamma}A_\alpha.\eqno{(4')}$$
Replace $2C_{\alpha\beta}-C_{\beta\alpha}$ and $2C_{\alpha\gamma}-C_{\gamma\alpha}$
in $(3'')$ by their expressions obtained from (4) and $(4')$, respectively. We obtain
$$C_{\beta\gamma}+C_{\gamma\beta}=C_{\beta\beta}A_\alpha(A_\beta^{-1}A_\gamma)A_\alpha^{-1}+
C_{\gamma\gamma}A_\alpha(A_\gamma^{-1}A_\beta)A_\alpha^{-1}.\eqno{(5)}$$

Take a point in $U$. In equation (5) written at this point, $\alpha$, $\beta$ and $\gamma$
may be arbitrary vectors. Fix $\beta$ and $\gamma$, and let $\alpha$
vary. Then $A_\alpha$ runs over all quaternions, since $f$ is a diffeomorphism.  
We need the following 

\begin{lemma}
Let $a$ and $b$ be arbitrary quaternions, and $x$ an imaginary quaternion (i.e. not 
a real number). Suppose that the number $ay+by^{-1}$ stays the same for all
quaternions $y$ obtained from $x$ by inner conjugations (i.e. $y=qxq^{-1}$ for
an arbitrary quaternion $q$). Then $b=a|x|^2$.
\end{lemma}

\proof 
Denote by $x_0$ the real part of $x$. Then it is easy to see that 
$$ay+by^{-1}=\left(a-\frac b{|x|^2}\right)y+\frac{2bx_0}{|x|^2}.$$
The second term in the right-hand side is independent of $y$. Therefore, 
the first term should be also independent of $y$ which is possible only if 
the coefficient $a-b/|x|^2$ vanishes. $\Box$

We can apply this lemma to equation (5) where we put $a=C_{\beta\beta}$, 
$b=C_{\gamma\gamma}$ and $x=A_\beta^{-1}A_\gamma$. If $\beta$ is not 
parallel to $\gamma$, then $x$ is imaginary. Thus by the lemma we have 
$$\frac{C_{\beta\beta}}{|A_\beta|^2}=\frac{C_{\gamma\gamma}}{|A_\gamma|^2}.$$
It follows that the left-hand side is independent of $\beta$. Hence
$$C_{\beta\beta}=C|A_\beta|^2\eqno{(6)}$$
where $C$ is a quaternion which does not depend on the vector $\beta$,
but it may depend on a point --- so $C$ is a quaternion-valued function
on $U$. 

Plug in (6) to (5):
$$C_{\beta\gamma}+C_{\gamma\beta}=C(\overline{A_\beta}A_\gamma+
\overline{A_\gamma}A_\beta).\eqno{(6')}$$
Comparing equations $(6')$ (with $\beta$ changed to $\alpha$ and $\gamma$ changed 
to $\beta$) and (4) we can conclude that 
$$C_{\alpha\beta}=\frac 13 C(\overline{A_\alpha}A_\beta+\overline{A_\beta}A_\alpha+
A_\alpha\overline{A_\beta}).\eqno{(6'')}$$
This is the {\em first integrability condition}. 

\section{Examples}

Let us now take a closer look on examples 1 and 2 form Section 1.

\paragraph{Example 1.}
Consider the left quaternionic Hopf fibration $\pi:\RP^7\to S^4$. 
A point in the preimage is represented by a pair of quaternions 
$(y,z)$ up to multiplication of both $y$ and $z$ by a common real factor. 
The image of this point under $\pi$ is an element of the left 
projective quaternionic line with homogeneous coordinates $y$ and $z$.
We can take $y^{-1}z$ as an affine conformal quaternionic coordinate in the 
image. 

Suppose that $f$ is obtained as the composition of $\pi$ with some 
projective embedding of (a part of) $\R^4$ to $\RP^7$. This projective 
embedding sends a point $x\in\R^4$ to a point in $\RP^7$ with coordinates
$$y=L(x),\quad z=M(x),$$
where $L$ and $M$ are some affine maps from $\R^4$ to $\R^4$. 
Thus in the given coordinates $f$ looks as follows: 
$$f:x\mapsto L(x)^{-1}M(x)\eqno{(E_1)}$$ 
where the multiplication and the inverse are in the sense of quaternions. 

Denote by $\vec L$ and $\vec M$ the linear parts of $L$ and $M$, respectively.
This means that $\vec L$ (resp., $\vec M$) is a linear operator such that 
$L$ (resp., $M$) is a composition of $\vec L$ (resp., $\vec M$) and a translation.
In other words, $\vec L$ and $\vec M$ are differentials of $L$ and $M$ at any 
point.

Compute $A_\alpha=\d_\alpha f$ for our map $f$:
$$A_\alpha(x)=-L^{-1}(x)\vec L(\alpha)L^{-1}(x)M(x)+L^{-1}(x)\vec M(\alpha).$$
Differentiating once again along $\alpha$ we obtain:
$$B_\alpha=-2L^{-1}(x)\vec L(\alpha).$$
We can see from here that $\d_\alpha B_\beta=\frac 12 B_\alpha B_\beta$, so
in this case $C=0$.

\paragraph{Case $C=0$.}
Let us work out the case when $C$ vanishes everywhere in $U$. In this case 
the first integrability condition $(6'')$ reads as follows:
$$\d_\alpha B_\beta=\frac 12 B_\alpha B_\beta.$$
Consider this equation together with equation (2). This is a system of partial 
differential equations of first order which expresses all first derivatives 
of $A$ and $B$ through the values of $A$ and $B$. Denote this system by $\mathbb{S}$. Fix any point 
$x\in U$ and initial values of 1-forms $A$ and $B$ at $x$. Thus for any vector
$\alpha$ at $x$ we know $A_\alpha$ and $B_\alpha$. Note that if a solution of 
$\mathbb{S}$ 
with given initial values exists, then it is unique. Indeed, all higher derivatives of $A$
and $B$ can be expressed through the initial values by means of the system $\mathbb{S}$. 
The existence of a solution follows from example 1. In this example, at the 
given point $x$ linear maps $A:\alpha\mapsto A_\alpha$ and $B:\alpha\mapsto B_\alpha$ can 
be arbitrary. 

We have just proved the following

\begin{prop}
\label{C=0}
If $C=0$ identically in $U$, then the general solution of equations (2) and 
$(6'')$ is given by $(E_1)$.
\end{prop}

\paragraph{Example 2.} 
Consider a Euclidean sphere $S$ in $\R^5$. 
We can choose this sphere to be centered at the origin and to have radius 1.
Introduce a conformal coordinate system on $S$ by means of the 
stereographic projection $j$ to the equatorial hyperplane 
with the center at the Noth pole.
To this end we need to choose an orthogonal splitting 
$\R^5=\R^4\times\R$ the latter factor being the line between the North and 
South poles. Thus any point in $\R^5$ is represented as $(y,z)$ where $y\in\R^4$
and $z\in\R$. The point on $S$ corresponding to $(y,0)$ by means of the 
stereographic projection is
$$j^{-1}(y)=\left(\frac{2y}{1+|y|^2},\frac{|y|^2-1}{1+|y|^2}\right).$$ 

Suppose that $f$ is the central projection from a horizontal hyperplane 
$z=z_1$ to $S$ with the center $(0,z_0)$. Under this projection, a point 
$(x,z_1)$ corresponds to a point $j^{-1}(y)$ on $S$ if and only if 
$$x=\frac{2y(z_0-z_1)}{(z_0+1)+(z_0-1)|y|^2}.$$

We can now rescale $y$ and choose specific values of $z_0$ and $z_1$ so that
$$x=\frac{y}{1+|y|^2}.\eqno{(E_2.1)}$$
This can be considered as an implicit equation defining $f:x\mapsto y$. 
Recall that $f$ was defined as a classical projection. 
Another interpretation of $f$ is that it establishes a correspondence between 
the Klein and the Poincar\'e models of the hyperbolic (Lobachevsky) geometry.  

Multiply both sides of equation $(E_2.1)$ by the denominator to get $y=(1+|y|^2)x$
and then differentiate both sides along $\alpha$: 
$$A_\alpha=\frac{2\<y,A_\alpha\>y}{1+|y|^2}+(1+|y|^2)\alpha.\eqno{(E_2.2)}$$
Here $\<\cdot,\cdot\>$ denotes the Euclidean inner product corresponding to 
the given quaternionic structure. Namely, for 2 vectors $\xi$ and $\eta$ we 
have $\<\xi,\eta\>=\Re(\xi\bar\eta)$ where in the right-hand side $\xi$ and 
$\eta$ are considered as quaternions. 
To express $A_\alpha$ it terms of $y$ only we need to find an expression for
$\<y,A_\alpha\>$. For that, take the inner product of both sides of equation 
$(E_2.2)$ with $y$. This leads to a linear equation for $\<y,A_\alpha\>$. Solving
this equation, we obtain
$$\<y,A_\alpha\>=\frac{(1+|y|^2)^2}{1-|y|^2}\<y,\alpha\>.\eqno{(E_2.3)}$$
Substitute $(E_2.3)$ into $(E_2.2)$. We now have a formula for $A_\alpha$
in terms of $y$: 
$$A_\alpha=(1+|y|^2)\left(\frac{2\<y,\alpha\>y}{1-|y|^2}+\alpha\right).\eqno{(E_2.4)}$$

Differentiating equation $(E_2.4)$ along $\alpha$ we can see that 
$$B_\alpha=2\left(2\frac{1+|y|^2}{1-|y|^2}\<y,\alpha\>+
\frac{y\overline{A_\alpha}}{1-|y|^2}\right).$$
Differentiate this equation once again along $\alpha$. Then we find that 
$$\d_\alpha B_\alpha-\frac 12 B_\alpha^2=\frac{6|A_\alpha|^2}{(1-|y|^2)^2}.$$
From equation (6) it follows that 
$$C=\frac 6{(1-|y|^2)^2}$$
in this example. In particular, $C$ may be nonzero. Note that in this example
$C$ is real everywhere. 

\section{Second Integrability Condition}

The case $C=0$ has been already worked out in the previous section. Now
assume that $C\ne 0$ at a given point. A priori, $C$ is a quaternion-valued
function. But now we are going to prove the following

\begin{lemma}
The function $C$ is real-valued.  
\end{lemma}

This lemma does not follow form equation (1). We need to use directly the fact 
that $f$ takes lines to circles. 

\proof
Take an arbitrary point $x_0\in U$ and a vector $\alpha$ at $x_0$. The line $l$
passing through $x_0$ in the direction of $\alpha$, is mapped to a circle under $f$.
Let $t\mapsto x=x(t)$ be the affine parameterization of $l$ such that $x(0)=x_0$ and 
$\dot x(0)=\alpha$. 

First assume that $B_\alpha$ is real at $x_0$. Then $f(l)$ must be a line.
Therefore, $\Im(B_\alpha)$ must vanish everywhere on $l$. It follows that 
$C$ is also real in this case.

Suppose now that $\Im(B_\alpha)\ne 0$. 
By Proposition 5.4 from \cite{Tim1} the center of the circle $f(l)$ is 
$$f(x)-(\Im B_\alpha(x))^{-1}A_\alpha(x).$$
Hence this expression does not depend on the choice of $x\in l$. Differentiate it 
by $t$:
$$A_\alpha+(\Im B_\alpha)^{-1}\left(\frac 12 \Im(B_\alpha^2)+\Im(C)|A_\alpha|^2\right)(\Im B_\alpha)^{-1}
A_\alpha-(\Im B_\alpha)^{-1}B_\alpha A_\alpha=0.$$
It follows that $\Im(C)=0$. $\Box$

The first integrability condition is a partial differential equation on $B$. 
Let us try to write down the integrability condition for this equation. 

We will use the language of quaternion-valued differential forms. By definition,
the algebra of quaternion-valued differential forms is obtained from the 
algebra of ordinary real-valued differential forms by taking the tensor product
with the quaternions over reals. In particular, this definition says how to 
multiply quaternion-valued forms. As was already mentioned, $A$ and $B$ are
thought of as quaternion-valued 1-forms. 

From the first integrability condition it follows that
$$dB=\frac 12 B\wedge B+\frac 13 C(A\wedge\bar A).\eqno{(6''')}$$
Take the differential of the both parts of this equation, using
relations $d^2B=0$ and $dA=0$:
$$dC\wedge A\wedge\bar A=\frac C2(B\wedge A\wedge\bar A-A\wedge\bar A\wedge B).\eqno{(7)}$$

\begin{prop}
\label{unq}
Suppose that the values of the 1-forms $A$ and $B$ and the function $C$ are given 
at some point, and the map $\alpha\mapsto A_\alpha$ is one-to-one at this point. 
If an analytic solution $f$ of equations (2) with these initial values of 
$A$, $B$ and $C$ exists, then it is unique. 
\end{prop}

\proof 
Since $A$ is non-degenerate, the operator of the right wedge multiplication by 
$A\wedge\bar A$ is invertible. This can be verified by a simple direct computation.
Therefore, equation (7) expresses all first derivatives of $C$ through $A$ and $B$.
Consider this equation together with equations (2) and $(6'')$. We obtain a 
system of partial differential equations which expresses all first derivatives 
of $A$, $B$ and $C$ in terms of the values of $A$, $B$ and $C$ only. Therefore, by successive 
applications of this system we can express all higher derivatives of $A$, $B$ and $C$
at the given point through the initial values at this point. The proposition
now follows. $\Box$ 

Note also that $A$, $B$ and $C$ determine the 3-jet of $f$. Namely, from 
equations (2) and $(6'')$ it follows that 
$$f=f(x_0)+A_{x-x_0}+\frac{B_{x-x_0}A_{x-x_0}}2+
\frac{(\frac 32B_{x-x_0}^2+C|A_{x-x_0}|^2)A_{x-x_0}}6+\dots$$
In particular, Proposition \ref{unq} can be reformulated as follows: 

\begin{prop}
\label{unq2}
If there exists a diffeomorphism $f:U\subseteq\R^4\to V\subseteq\R^4$ with a given 3-jet at 
some point $x_0\in U$ such that the image of each line segment lying in $U$ is an arc of a 
circle or a line segment in $V$, then such diffeomorphism is unique. 
\end{prop}

\section{Admissible 3-jets}

From now on we identify the space of preimage with the space of image and 
assume that the distinguished point is the origin, $f(0)=0$ and $A_\alpha=\alpha$ for all
vectors $\alpha$ at 0. This may be achieved by a linear change of variables in the
preimage composed with a translation in the image. 
Thus if $x=x^0+ix^1+jx^2+kx^3$ denotes the natural quaternionic coordinate, 
then $A=dx=dx^0+idx^1+jdx^2+kdx^3$. 

Equation (7) imposes a restriction on $B$ due to the fact that $C$ is real. 
Let us find all admissible values of $B$ at 0 by solving (7) as a linear
system on $B$ and $dC/C$. We have 
$A\wedge\bar A=dx\wedge\overline{dx}=i\omega_1+j\omega_2+k\omega_3$
where $\omega_1$, $\omega_2$ and $\omega_3$ are real 2-forms given in 
coordinates by
$$
\begin{array}{c}
\omega_1=2(dx^1\wedge dx^0+dx^3\wedge dx^2),\\
\omega_2=2(dx^2\wedge dx^0+dx^1\wedge dx^3),\\
\omega_3=2(dx^3\wedge dx^0+dx^2\wedge dx^1).
\end{array}
$$
Denote by $\gamma$ the real 1-form $2dC/C$ and suppose that 
$B=B^0+iB^1+jB^2+kB^3$ where $B^0$, $B^1$, $B^2$ and $B^3$ are real 1-forms. 
Equation (7) can be now rewritten as the following system: 
$$\begin{array}{c}
2(B^2\wedge\omega_3-B^3\wedge\omega_2)=\gamma\wedge\omega_1,\\
2(B^3\wedge\omega_1-B^1\wedge\omega_3)=\gamma\wedge\omega_2,\\
2(B^1\wedge\omega_2-B^2\wedge\omega_1)=\gamma\wedge\omega_3.
\end{array}\eqno{(8)}$$
Denote by $B^\mu_\nu$ and $\gamma_\nu$ ($\nu=0,1,2,3$) the components of 
the 1-forms $B^\mu$ and $\gamma$ so that 
$$B^\mu=\sum_{\nu=0}^3 B^\mu_\nu dx^\nu,\quad \gamma=\sum_{\nu=0}^3 \gamma_\nu dx^\nu.$$
System (8) yields the following equations: 
$$\begin{array}{c}
2(-B^2_3+B^3_2)=-\gamma_1,\ 2(-B^2_2-B^3_3)=\gamma_0,\ 2(B^2_1-B^3_0)=-\gamma_3,\
2(B^2_0+B^3_1)=\gamma_2,\\
2(-B^3_1+B^1_3)=-\gamma_2,\ 2(B^3_0+B^1_2)=\gamma_3,\ 2(-B^3_3-B^1_1)=\gamma_0,\
2(B^3_2-B^1_0)=-\gamma_1,\\
2(-B^1_2+B^2_1)=-\gamma_3,\ 2(B^1_3-B^2_0)=-\gamma_2,\ 2(B^1_0+B^2_3)=\gamma_1,\
2(-B^1_1-B^2_2)=\gamma_0.
\end{array}$$
Exclude $\gamma_\nu$ from these equations: 
$$\begin{array}{c}
B^1_1=B^2_2=B^3_3,\\
B^1_0=B^2_3=-B^3_2,\\
B^2_0=B^3_1=-B^1_3,\\
B^3_0=B^1_2=-B^2_1.
\end{array}$$
The relations displayed above mean that $B_\alpha=p(\alpha)+\alpha q$
where $p$ is a real-valued 1-form and $q$ is a quaternion. Hence the 3-jet 
at 0 of $f$ should be as follows: 
$$x+\frac{(p(x)+xq)x}2+\frac{(\frac 32(p(x)+xq)^2+C|x|^2)x}6.\eqno{(9)}$$

\begin{prop}
\label{ex}
For an arbitrary real-valued linear function $p$ on $\R^4$ and an 
arbitrary quaternion $q$,
a local diffeomorphism with 3-jet (9) exists. 
Moreover, it can be chosen to be one of the classical projections.
\end{prop}

\proof 
First consider a map given by the implicit equation 
$$x=\frac y{1-\frac C6|y|^2}.$$
This is clearly a map from example 2.
It has the 3-jet $x+\frac C6|x|^2x$. This is a partial case of (9) where 
$p=q=0$. To achieve other given values of $p$ and $q$, compose this map with the following 
M\"obius transformation: 
$$y\mapsto 2q^{-1}\left(1-\frac{\frac 12 qy}{1-\frac 12p(y)}\right)^{-1}-2q^{-1}.$$
It is easy to see that the composition has 3-jet (9). $\Box$

\paragraph{Concluding remarks.} We should now make things add up in the proof of 
theorem \ref{main}.  

We have always assumed that $\d_\alpha A_\alpha=B_\alpha A_\alpha$ the multiplication
by $B_\alpha$ being from the left. The case of the right multiplication is completely 
analogous. The only difference is that for $C=0$ we would have the right quaternionic
Hopf fibration instead of the left one. In general there is an open subset $U'$ of $U$ 
such that the multiplication by $B_\alpha$ is either form the left everywhere on $U'$ or from the 
right everywhere on $U'$. Therefore theorem \ref{main} holds on $U'$. By the uniqueness
theorem for analytic functions it holds then on $U$. 

If $C=0$ everywhere on $U$, then by Proposition \ref{C=0} we have a quaternionic Hopf
fibration. Otherwise $C$ nowhere vanishes on some open subset of $U$. Then by Propositions
\ref{unq2} and \ref{ex} the map $f$ is a classical projection on this subset. By the 
uniqueness theorem $f$ is a classical projection on $U$.

\bigskip
\bigskip

{\sc Department of Mathematics, University of Toronto
100 St. George street, Toronto ON M5S3G3 Canada}

{\em E-mail address:} {\sf vtimorin@math.toronto.edu}

\end{document}